\documentclass[twoside,12pt]{article}

\usepackage{amsmath}
\usepackage{amssymb}
\usepackage{amscd}
\usepackage{amsthm}

\numberwithin{equation}{section}

\theoremstyle{plain}
\newtheorem{theorem}{Theorem}[section]

\theoremstyle{remark}
\newtheorem{remark}[equation]{Remark}
\theoremstyle{definition}

\newcommand{\R}{\mathbb R}
\newcommand{\C}{\mathbb C}

\newcommand{\Q}{\mathbb Q}

\def\ra{\rightarrow}

\def\e{\emph}
\def\i{\infty}

\def\b{\begin}

\begin{document}

\title{
{Rigidity of   fiber-preserving quasisymmetric maps }}
\author{Enrico Le Donne  and    Xiangdong Xie
}
\date{January 7, 2015}

\maketitle

\begin{abstract}
We show that fiber-preserving quasisymmetric maps are biLipschitz. As an application, we show that  quasisymmetric maps on 
 Carnot groups 
  with reducible first stratum    are biLipschitz.


\end{abstract}




 {\small {\bf{Mathematics Subject
Classification (2010).}} 30L10,  22E25,     53C17. 






\setcounter{section}{0} 
  \setcounter{subsection}{0}

\section{Introduction}\label{s0}

In this paper we study the rigidity property of quasisymmetric maps that preserve a foliation. 
 We show that, under quite general conditions, such
   quasisymmetric maps are biLipschitz. 
  We then give an application
     to  quasisymmetric maps between Carnot groups.

Quasisymmetric maps that preserve a foliation arise  when one studies the rigidity  property of quasiisometries between negatively curved solvable Lie groups, see  \cite{Pi},  \cite{SX},  
\cite{X1}, \cite{X4}.  It is well known that a negatively curved space has an
 ideal  boundary, and quasiisometries between negatively curved spaces correspond to quasisymmetric maps between the ideal boundaries.   Rigidity properties of the quasiisometries correspond to the rigidity properties of the  quasisymmetric maps. The ideal boundary of a negatively curved solvable Lie group is a nilpotent Lie group 
 equipped with a homogeneous distance.  Recent results suggest that very often the quasisymmetric maps 
on such nilpotent Lie groups must be biLipschitz \cite{P}, \cite{CO},  \cite{CR}, \cite{SX}, 
      \cite{X1},  \cite{X2},   \cite{X4}, \cite{X5}.  

There are usually two steps in the   proof of the above mentioned rigidity property of the quasisymmetric maps. 
 The first step is to show that the quasisymmetric map preserves a certain foliation. The second step is to show 
 that  a quasisymmetric map must be biLipschitz if it preserves a foliation.  The main purpose of this paper
 is   to take a  closer look at the second step  in the context of quasisymmetric maps between general metric spaces.   We hope that this will be useful in the eventual  (hopefully) complete solution of the rigidity of quasiisometries between negatively curved solvable Lie groups.

One often has to deal with quasimetrics instead of metrics while studying the ideal boundary of negatively curved spaces. For this reason we state our main result 
 for quasimetric spaces.   Recall that a function 
 $d: X\times X\ra [0, \infty)$ is a {\it{quasimetric}}  on a   set $X$   if \newline
 (1) $d$ is symmetric, that is,  $d(x_1, x_2)=d(x_2, x_1)$  for all $x_1, x_2\in X$;\newline
 (2) $d(x_1, x_2)=0 $ if and only if $x_1=x_2$;\newline
  (3)    there exists some constant $M\ge 1$ such that 
  $d(x_1, x_3)\le M\cdot (d(x_1, x_2)+d(x_2, x_3))$ for all $x_1, x_2, x_3\in X$.  \newline
 We  call $M$   a  quasimetric constant of $d$.
  A subset $A$ of a quasimetric space $X$ is called {\it{closed}} if for any 
 $x\in X$  and any sequence $\{a_j\}\subset A$, the condition $d(a_j, x)\ra 0$ implies $x\in A$.

To state our main result,   we introduce the  notion of ``fibered
quasimetric space''.
  To do that we first  specify the concept of
parallelism for sets.
We say that two closed subsets  $U$ and $V$  of a quasimetric space $X$ are {\em parallel} if 
 there is a constant $a>0$ such that 
  $d(u, V)=d(v, U )=a$ for every 
$u\in U $ and every $v\in V $. 
 Recall that $d(u,V):=\inf\{d(u,v):v\in V\}$.
  It is easy to check that in this case we have 
 $d(U, V)=HD(U, V)=a$, where 
$$d(U,V):=\inf\{d(u,v)| u\in U,\; v\in V\}, $$
   and 
$HD(U, V)$ denotes the  Hausdorff distance between $U$ and $V$: 
$$ HD(U,V):= \sup (\{ d(u,V):u\in U\}\cup \{d(v,U):v\in V\}).$$

The following definition   provides the setting for our main theorem.


\b{Def}[Fibered   quasimetric space]
\label{fibered}
Let $\alpha>0$, $L\ge 1$. We say  that a quasimetric space $X$ is  an $(\alpha, L)$-{\em fibered quasimetric space}   if  
$X$ admits a 
cover $\mathcal U$ by  closed pairwise   disjoint
subsets, called {\em fibers},
   with   the following   properties:\newline
 \begin{equation}\label{fasfetugs}
 \text{
  \hspace{-3cm}
 Fibers are snow-flake equivalent to unbounded geodesic spaces:    }
 \end{equation}
    for each $U\in \mathcal U$, there exists an unbounded   geodesic space $(\widetilde{U}, d)$ 
 such that  $U$ is  $L$-biLipschitz to 
  $(\widetilde{U}, d^\alpha)$;  
 \begin{equation}\label{pfant}
  \hspace{-8.2cm}
 \text{Parallel fibers are not isolated:   }
 \end{equation}
for any $U\in \mathcal U$,  there exists a sequence  
$U_{i} \in \mathcal U$ such that $U_i$  and $U$  are distinct,  parallel, and 
$HD(U_i,U)\to 0$ as $i\to \infty$;
\newline
 \begin{equation}\label{fhpd}
  \hspace{-8.5cm}
 \text{
 Fibers have positive distance:   }
 \end{equation}
for any two distinct   $U,V\in\mathcal U$, 
  $d(U,V)>0$;\newline 
 \begin{equation}\label{npfd}
 \hspace{-8.8cm}
 \text{Non-parallel fibers diverge:   }
 \end{equation}
if  $U,V\in \mathcal U$ are not parallel, then  
$  HD(U,V)=\infty$.
\end{Def}

Given $K\ge 1$ and $C>0$, we say that a bijection $F:X\ra Y$ between two
     quasimetric spaces is  a $(K, C)$-\e{quasi-similarity}  if
\[
   \frac{C}{K}\cdot d(x_1,  x_2)\le d(F(x_1), F(x_2))\le CK\cdot d(x_1,  x_2), \qquad \forall x_1, x_2 \in X.
\]
Clearly a map is a quasi-similarity if and only if it is biLipschitz.  
 The point here is that often there is control on $K$  but not on $C$.  In this case, the notion of quasi-similarity provides more information about the distortion.

The following is our main result and states that any  
fiber-preserving quasisymmetric map between fibered   quasimetric spaces is a 
quasi-similarity (and hence biLipschitz). See Section~\ref{maps} for the definition of quasisymmetric map.

\b{theorem}\label{main1}
Let  $X, Y$ be  $(\alpha, L)$-fibered  quasimetric spaces for some $\alpha>0 $ 
and $L\ge 1$.
   Suppose $F: X\ra Y$ is an  $\eta$-quasisymmetric map that sends fibers of $X$ homeomorphically onto fibers of $Y$.  Then $F$ is a $(K, C)$-quasi-similarity,
 where $K$ depends only on $\eta$, $\alpha$,  $L$  and the quasimetric constants of $X$, $Y$.

\end{theorem}

  We remark that 
 for the validity of  Theorem \ref{main1} the condition $d(u, V)=d(v, U)=a$  for  parallel fibers can not be  replaced by the weaker  condition $HD(U, V)<\infty$. See the end of Section  
 \ref{fiber}  for an example.


Theorem~\ref{main1} in particular implies that all quasiconformal maps of the subRiemannian Heisenberg groups that send vertical lines to vertical lines must be biLipschitz, see  Proposition~\ref{hei}.     See \cite{T} and \cite{BHT}  for the construction of such maps. 
 We remark that there exist biLipschitz maps of  Heisenberg groups that
 send vertical lines to curves that are not vertical lines, see  \cite{X3}.

Our main application  of Theorem~\ref{main1}
   is to quasisymmetric maps on Carnot groups    with reducible first stratum.

 Let  $G$ be a  Carnot group.
 See Section~\ref{basics} for a brief introduction to Carnot groups.
Let    $V_1\oplus \cdots \oplus V_r$ be the stratification of the Lie algebra
$\mathfrak g  $
of $G$.
  For each $t>0$, the   standard  dilation  $\lambda_t: \mathfrak g \ra \mathfrak g $ is defined 
 by $\lambda_t(v)=t^i v$ for $v\in V_i$.
 A  Lie algebra    isomorphism $A:  \mathfrak g \ra \mathfrak g $ is  called a {\it{strata-preserving}} automorphism
   if 
   $A(V_i)=V_i$, for $i=1, \ldots, r$.
   Let $\text{Aut}_*(\mathfrak g )$ be the group of strata-preserving automorphisms of $\mathfrak g $. We say $V_1$ is 
{\it reducible} (or the first stratum of $\mathfrak g $ is reducible) if  there is a non-trivial proper
 linear subspace $W_1\subset V_1$ such that $A(W_1)=W_1$ for every 
 $A\in \text{Aut}_*(\mathfrak g )$. 

 
The main consequence of Theorem~\ref{main1} is the following theorem.
Recall that a Carnot group is a  metric space equipped with a particular Carnot-Carath\'eodory distance, see Section~\ref{basics}.

\b{theorem}\label{main2}
Let $G$ be  a Carnot group with reducible first stratum.
Then every quasisymmetric map $F: G\ra G$ is a quasi-similarity, quantitatively.
\end{theorem}
Theorem~\ref{main2} is quantitatively in the sense that 
for any Carnot group $G$ and for any 
$\eta\in {\rm Homeo}([0,\infty))$ there exists a constant $K$ depending 
  only   on $G$ and $\eta$
such that every
$\eta$-quasisymmetric map
$F: G\ra G$ is a 
$(K, C)$-\e{quasi-similarity} for some $C$.


We remark that in \cite{CO}, Cowling and Ottazzi proved that all smooth quasiconformal maps on rigid Carnot groups are affine maps (an affine map is the composition of a left translation and a
 strata-preserving
 automorphism) and so in particular are biLipschitz. On the other hand, 
   it is shown in \cite{X5} that  if $N$ is a non-rigid Carnot group other than the Euclidean groups  and some very particular type of product groups,  then every quasisymmetric map  of $N$
 is  biLipschitz. 
  We note that Carnot groups with reducible first stratum include both rigid and non-rigid groups. 
 In Theorem \ref{main2} there is no assumption on the regularity of the quasisymmetric map.

In Section~\ref{prelimi} we recall some basic definitions and facts that are used later. 
 In Section~\ref{fiber}   we prove   Theorem~\ref{main1}.  In Section~\ref{app}  we 
   prove Theorem~\ref{main2} and  also give  two other applications of Theorem~\ref{main1}:   quasiconformal maps of Heisenberg groups that send vertical lines to vertical lines, and   a new proof of a theorem of Dymarz  \cite{D}    on quasisymmetric maps of   ideal boundary  of 
  certain  amenable hyperbolic locally compact groups.

\noindent {\bf{Acknowledgment}}. {This work was initiated when the second-named author visited  the University of 
Jyvaskyla.  He  would like to thank  University of 
Jyvaskyla  for  financial support  and  excellent working conditions.
   He  also  acknowledges support from
 NSF grant  DMS-1265735
     and   Simons Foundation grant \#315130.  }


\section{Preliminaries}\label{prelimi}

In this section we collect definitions and results   needed later. 
  We first recall the 
definitions  of   quasiconformal  and   quasisymmetric maps  (Section~\ref{maps}). 
  Then we  review the 
basic definitions related to Carnot groups  in Section~\ref{basics}, 
   the BCH formula in Section~\ref{BCH}   
and    
  Pansu Differentiability Theorem in Section~\ref{pansud}.

  \subsection{Quasiconformal  and quasisymmetric   maps}\label{maps}


 A map $F: X \ra Y$ between two quasimetric spaces is continuous if for any $x\in X$ and any sequence $\{x_j\}\subset X$, the condition $\lim_{j\ra\infty} d(x, x_j)=0$ implies  $\lim_{j\ra\infty} d(F(x), F(x_j))=0$.

Let $F:  X\ra Y$ be a  bijection  between two   quasimetric spaces.
 For $x\in X$ and $t>0$, define
  $$H_F(x,t)=\frac{\sup\{d(F(x'), F(x))| d(x', x)\le t\}}{\inf\{d(F(x'), F(x))| d(x', x)\ge t\}}.$$
 The map   $F$  is called $\lambda$-{\em quasiconformal} if 
  both  $F$ and $F^{-1}$ are continuous and 
 $$\limsup_{t\ra 0} H_F(x,t)\le \lambda$$
   for all $x\in X$. 
    We say  $F$ is quasiconformal if it is $\lambda$-quasiconformal for some $\lambda\ge 1$.

Let $\eta: [0,\i)\ra [0,\i)$ be a homeomorphism.
    A   bijection 
$F:X\to Y$ between two  quasimetric spaces is
\e{$\eta$-quasisymmetric} if 
both  $F$ and $F^{-1}$ are continuous and 
for all distinct triples $x,y,z\in X$,
we have
\[
   \frac{d(F(x), F(y))}{d(F(x), F(z))}\le \eta\left(\frac{d(x,y)}{d(x,z)}\right).
\]
  If $F: X\ra Y$ is an $\eta$-quasisymmetry, then
  $F^{-1}:    Y\ra X$ is an $\eta_1$-quasisymmetry, where
$\eta_1(t)=(\eta^{-1}(t^{-1}))^{-1}$, see \cite[Theorem 6.3]{V}.
 A     bijection  between quasimetric spaces
  is said to be quasisymmetric if it is $\eta$-quasisymmetric for some $\eta$. 

We remark that 
     quasisymmetric  maps between general quasimetric spaces
  are quasiconformal. In the case of Carnot groups,  and more
  generally  Loewner  spaces,  a   quasiconformal homeomorphism  is locally  quasisymmetric, see \cite[Theorem 4.7]{HK}.

\subsection{Stratified Lie  algebras and Carnot groups}\label{basics}


A   \e{stratified Lie algebra} is a finite-dimensional Lie algebra  
$\mathfrak g$  over $\R$      together with  a direct sum   decomposition  
    $\mathfrak g=V_1\oplus V_2\oplus\cdots \oplus V_s$ 
  of      nontrivial  vector subspaces
 such that $[V_1,
V_i]=V_{i+1}$ for all $1\le i\le s$,
    where we set $V_{s+1}=\{0\}$.  The integer $s$ is called step or   
    degree of nilpotency of $\mathfrak g$. Every stratified Lie algebra
 $\mathfrak g=V_1\oplus V_2\oplus\cdots \oplus V_s$  admits a one-parameter
 family of automorphisms $\lambda_t: \mathfrak g \ra \mathfrak g$, $t\in (0, \i)$, where
 $\lambda_t(x)=t^i x$ for  $x\in V_i$.

  Let   $\mathfrak g=V_1\oplus V_2\oplus\cdots \oplus V_s$
    and $\mathfrak g'=V'_1\oplus V'_2\oplus\cdots \oplus V'_s$  be two  stratified Lie
    algebras.
  A Lie algebra  isomorphism
  $\phi: \mathfrak g\ra \mathfrak g'$
  is such that 
   $\phi(V_i)=V'_i$ for all $1\le i\le
  s$
  if and only if
  it commutes with $\lambda_t$ for
  all $t>0$; that is, if $\phi\circ \lambda_t=\lambda_t\circ
  \phi$.  
  We call such isomorphisms  {\it{strata preserving}}.

Let $G$ be a connected,  simply connected Lie group 
 whose  Lie algebra is stratified as
      ${\rm Lie}( G)=V_1\oplus \cdots \oplus V_s$.  
      The subspace $V_1$ defines
      a left-invariant distribution  $H G\subset TG$ on $G$.    We fix a left-invariant inner product on
          $HG$.
           An
      absolutely continuous curve $\gamma$ in $G$  whose velocity vector
       $\gamma'(t)$  is contained in  $H_{\gamma(t)} G$ for almost every $t$
        is called  a {\em horizontal curve}.
        Any horizontal curve has an associated length defined using the left-invariant inner product, i.e., by integrating the norm of its tangent vector.
  Since $V_1$ generates the whole Lie algebra,         any two points
  of $G$ can be  connected by horizontal curves. Let $p,q\in G$, the
  \e{Carnot-Carath\'eodory  metric} between $p$ and $q$ is denoted by  $d_c(p,q)$ and      is defined as
  the infimum of length of horizontal curves that join $p$ and $q$.
  We call {\em Carnot group} the data of $G$, its stratification, the inner product on its first stratum, and consequently the distance $d_c$.

  Since the inner product on   $HG$ is left invariant, the Carnot
  metric on $G$ is left invariant as well.  Different choices of inner
  product on $HG$ result in Carnot metrics that are biLipschitz
  equivalent.
    The Hausdorff dimension of $G$ with respect to  a  Carnot metric
    is  given by $\sum_{i=1}^s i\cdot \dim(V_i)$.

Recall that, for a simply connected nilpotent Lie group $G$ with Lie
algebra $\mathfrak g$, the exponential map
  $\exp: {\mathfrak g}\ra G$ is a diffeomorphism.  Under this identification the Lesbegue   measure on 
  $\mathfrak g$  is a  Haar measure on $G$.  
 Furthermore, the
  exponential map induces 
  a  one-to-one correspondence between
    Lie subalgebras of $\mathfrak g$   and
  connected Lie subgroups of $G$.

Let $G$ be a Carnot group with Lie algebra
      $\mathfrak g=V_1\oplus \cdots \oplus  V_s$.
         Since $\lambda_t: \mathfrak g\ra  \mathfrak g$ ($t>0$) is a Lie algebra
         automorphism  and $G$ is simply connected,  there is a
         unique  Lie group automorphism  $\Lambda_t: G\ra G$ whose
         differential  at the identity is $\lambda_t$.
       For each $t>0$,
      $\Lambda_t$ is a
             similarity with respect to the Carnot metric:  $d(\Lambda_t(p),
             \Lambda_t(q))=t\, d(p,q)$ for any two points $p, q\in
             G$.  
             For a  Lie group isomorphism
                $f: G\ra G'$ between two Carnot groups,
                the corresponding Lie algebra
                 isomorphism
                 $f_*: {\rm Lie} ( G)\ra {\rm Lie} ( G')$
                 is strata preserving 
                 if and only if
                  $f$ commutes with
                    $\Lambda_t$ for all  $t>0$; that is, if
                      $f\circ \Lambda_t=\Lambda_t\circ  f$.

\subsection{The  Baker-Campbell-Hausdorff  formula}\label{BCH}

Let  $G$ be a simply connected nilpotent Lie group with Lie algebra $\mathfrak g$.
  The exponential map  $\text{exp}:  {\mathfrak g}\ra G$ is a diffeomorphism.  
   One can then
     pull back the group operation from $G$ to get a group structure 
 on $\mathfrak g$.     This group structure can be described by the   Baker-Campbell-Hausdorff formula
 (BCH formula in short),  which expresses the product $X*Y$ ($X, Y\in {\mathfrak g}$)
    in terms of the iterated Lie brackets  of $X$ and $Y$. 
     The group operation in $G$ will be denoted by $\cdot$.  
   The pull-back group operation   $*$  on $\mathfrak g$ is defined as follows. 
      For $X,   Y\in \mathfrak g$,   define
   $$X*Y=\text{exp}^{-1}(\text{exp} X\cdot \text{exp} Y).$$
  Then the first a few terms of the BCH formula (\cite{CG},    page 11)  is given by:
$$ X*Y =X+Y+\frac{1}{2}[X,Y]+\frac{1}{12}[X,[X,Y]]
-\frac{1}{12}[Y, [X, Y]]+\cdots. $$

\subsection{Pansu differentiability  theorem }\label{pansud}

First the definition:

\b{Def}\label{pansu-d}
 Let $G$ and $G'$
  be two Carnot groups endowed with Carnot metrics,  and $U\subset G$, $U'\subset G'$ open subsets.
   A map $F: U\ra U'$ is \e{Pansu  differentiable}
    at $x\in U$  if there exists a strata preserving  homomorphism
     $L: G\ra G'$ such that
     $$\lim_{y\ra x}\frac{d(F(x)^{-1}*F(y),\, L(x^{-1}*y))}{d(x,y)}=0.$$
          In this case, the strata preserving homomorphism
     $L: G\ra G'$  is called the \e{Pansu  differential} of $F$ at $x$, and
     is denoted by   $dF(x)$.
\end{Def}

The following  result (except the terminology)   is due to  Pansu
[P].

\b{theorem}\label{pan}
  Let $G, G'$ be   Carnot groups, and $U\subset G$, $U'\subset G'$  open subsets.  
 Let $F: U\ra U'$ be a quasiconformal  map.
   Then $F$ is almost everywhere Pansu  differentiable. Furthermore, at a.e.
   $x\in   U$, the Pansu  differential $dF(x): G\ra G'$ is a strata preserving
   isomorphism.

\end{theorem}


\section{Fiber-preserving  quasisymmetric   maps  }\label{fiber}

In this section we prove  Theorem~\ref{main1}.

Let $X, Y$ be $(\alpha, L)$-fibered quasimetric spaces and $F: X\ra Y$ an $\eta$-quasisymmetric map that   sends   fibers  of $X$ onto fibers of $Y$. 
 We shall show that $F$ is a  $(K, C)$-quasi-similarity, where $K$ depends only on $\eta$, $\alpha$,   $L$  and the quasimetric constants of $X$, $Y$.

\b{Le}\label{parallel}
If two fibers  $U$  and $V$ in $X$  are parallel, then  
$F(U)$  and $F(V)$ are parallel fibers  in $Y$.
\end{Le}

\b{proof}
Suppose 
$F(U)$  and $F(V)$  are not   parallel.
 By Condition \eqref{npfd},   after possibly switching $U$ and $V$, 
   we may assume that there exist a sequence $x_i\in U$ such that 
 ${d(F(x_i), F(V))}\ra \infty$.  
  Since $U$  and $V$  are parallel, for each $i$,  there exists  $v_i\in V$
 such that $d(x_i, v_i)=d(U, V)$. 
   Since by Condition \eqref{fasfetugs}
   the set $U$ is path connected and unbounded,  
 there are points on $U$ at arbitrary distance from any point. Hence, for each $i$, there exists $u_i\in U$  
 such that $d(x_i, u_i)=d(x_i, v_i)$. 
 
Since        ${d(F(x_i), F(V))}\ra \infty$,    we  have 
$d(F(x_i), F(v_i))\ra \i$.  
  The quasisymmetry condition  and the fact 
$ d(x_i, u_i)=d(x_i, v_i)$
 imply  that 
 $d(F(x_i), F(u_i))\ra \i$.
 
 Because of Condition \eqref{pfant}, we can take a  fiber  $U'\not=F(U)$  parallel to 
 $F(U)$.  
Again by Condition \eqref{fasfetugs}
 choose $y_i\in U'$ so that
 $d(F(x_i),  y_i)=d(F(U),U')$.  Then 
 $\frac{d(F(x_i),  y_i)}{d(F(x_i),  F(u_i))}\ra 0$ as $i\ra \i$.  Now the quasisymmetry condition for $F^{-1}$ implies  
 $$\frac{d(x_i, F^{-1}(y_i))}{d(U,  V)}=
\frac{d(x_i, F^{-1}(y_i))}{d(x_i, u_i)}\ra 0 \; \;\text{as}\; \; i\ra \infty.$$
  It follows that ${d(x_i, F^{-1}(y_i))}\ra 0$. Since  $x_i\in U$ and
$F^{-1}(y_i)\in F^{-1}(U')$, we have 
 $d(U,  F^{-1}(U'))=0$, contradicting    Condition \eqref{fhpd}  
     and the fact   $U'\not=F(U)$. 
\end{proof}

The next two results  are similar to  Lemma 15.3  and Corollary 15.4 in \cite{PW}.
    The proofs  are    modifications of their arguments.

\b{Le}\label{l1}
 There exists  $K_1\ge 1$   depending only on $\eta$,  $\alpha$ and $L$ so that
 for any two parallel fibers $U_1$, $U_2$  in $X$ 
  and any $p, q\in U_1$ 
  satisfying  $d(p,q)\ge L\cdot  d(U_1, U_2)$  and  $d(F(p),  F(q))\ge L\cdot  d(F(U_1),   F(U_2))$  
 we have 
 \b{equation}\label{e0}\frac{1}{K_1}\cdot \frac{d(F(U_1), F(U_2))}{d(U_1, U_2)}\le  \frac{d(F(p), F(q))}{d(p, q)}
\le  K_1  \cdot\frac{d(F(U_1), F(U_2))}{d(U_1, U_2)}. 
\end{equation}

\end{Le}
\b{proof}
By  Condition \eqref{fasfetugs},     there exists  a  geodesic space $(\widetilde{X}_1, d_1)$   and 
 an  $L$-biLipschitz map $f_1: (U_1, d)\ra (\widetilde{X}_1, d^\alpha_1)$.  
  The  assumption  $d(p,q)\ge L\cdot  d(U_1, U_2)$   implies
   $$d_1(f_1(p), f_1(q))\ge  \left( d(U_1, U_2)\right)^{1/{\alpha}}.$$
 There is some integer $k\ge 2$ such that 
$$(k-1)\cdot   \left( d(U_1, U_2)\right)^{1/{\alpha}}\le  d_1(f_1(p), f_1(q)) < k\cdot  \left( d(U_1, U_2)\right)^{1/{\alpha}}.$$
  The  biLipschitz property  of $f_1$ implies
\b{equation}\label{e1}
\frac{1}{L} (k-1)^\alpha \cdot d(U_1, U_2)\le d(p,q)< L k^\alpha \cdot d(U_1, U_2).
\end{equation}
 Since $(\widetilde{X}_1, d_1)$ is a geodesic space, 
 there are points $p=p_0, p_1, \ldots, p_k=q$  in $U_1$ such that
 $\frac{1}{2}d(U_1, U_2)^{\frac{1}{\alpha}}\le d_1(f_1(p_i), f_1(p_{i+1}))\le 
d(U_1, U_2)^{\frac{1}{\alpha}}$.  
We  have 
 $$\frac{1}{L 2^\alpha} \cdot  d(U_1, U_2)\le d(p_i, p_{i+1})\le L \cdot  d(U_1, U_2).$$
  Since $U_1$ and $U_2$ are parallel,  Lemma~\ref{parallel} implies that 
 $F(U_1)$ and $F(U_2)$ are also parallel.  
   Let $q_i\in U_2$ be a point  such that $d(F(p_i), F(q_i))=d(F(U_1), F(U_2))$.
  We have 
$$d(p_i, p_{i+1})\le  L \cdot  d(U_1, U_2)  \le L \cdot  d(p_i,   q_i).$$
 Since  $F$ is $\eta$-quasisymmetric,    we have 
$$d(F(p_i), F(p_{i+1}))\le  \eta(L)\cdot d(F(p_i), F(q_i))=
 \eta(L) \cdot d(F(U_1), F(U_2)).$$
By  Condition \eqref{fasfetugs}  again,  
  there exists   some geodesic space $(\widetilde Y_1, \rho_1)$ and an  $L$-biLipschitz map 
 $g_1: (F(U_1), d)\ra (\widetilde Y_1, \rho^\alpha_1)$. 
  It follows that 
 $$\rho_1(g_1\circ F(p_i), \;  g_1\circ F(p_{i+1}))\le (L\eta(L))^{\frac{1}{\alpha}}\cdot  d(F(U_1), F(U_2))^{\frac{1}{\alpha}}.$$
 The triangle inequality   for the metric space $(\widetilde Y_1, \rho_1)$
implies 
  $$\rho_1(g_1\circ F(p), \; g_1\circ F(q))\le k  (L\eta(L))^{\frac{1}{\alpha}}\cdot  d(F(U_1), F(U_2))^{\frac{1}{\alpha}}.$$
  Now the biLipschitz property of   $g_1$ implies
  \b{equation}\label{e2}
d(F(p), F(q))\le L \cdot \rho^\alpha_1(g_1\circ F(p), g_1\circ F(q))\le L^2 \eta(L)k^\alpha
    \cdot  d(F(U_1), F(U_2)).
\end{equation}
Now the   second inequality in (\ref{e0})   follows from (\ref{e1})  and (\ref{e2}).
Finally, we notice that the second inequality for $F^{-1}$ is equivalent to the first inequality for $F$. 
\end{proof}

\b{Le}\label{l2}
  The following holds 
 for any  two  distinct parallel fibers $U_1$, $U_2$ 
  and any  two distinct $p, q\in U_1$:
 \b{equation}\label{e4}
\frac{1}{K_1^3}\cdot \frac{d(F(U_1), F(U_2))}{d(U_1, U_2)}\le  \frac{d(F(p), F(q))}{d(p, q)}
\le  K_1^3  \cdot    \frac{d(F(U_1), F(U_2))}{d(U_1, U_2)},
\end{equation}
where $K_1$ is the constant in Lemma~\ref{l1}.

\end{Le}

\b{proof}  Let $U_1$ and $U_2$ be two distinct parallel fibers,  and $p, q\in U_1$ be distinct.
 Pick  two points $p_0, q_0\in U_1$  that satisfies
 $d(p_0, q_0)> L\cdot d(U_1, U_2)$  and  $d(F(p_0), F(q_0))>L\cdot   d(F(U_1), F(U_2))$.  
By Lemma~\ref{l1} we have 
\b{equation}\label{e5}
\frac{1}{K_1}\cdot \frac{d(F(U_1), F(U_2))}{d(U_1, U_2)}\le  \frac{d(F(p_0), F(q_0))}{d(p_0, q_0)}
\le  K_1\cdot    \frac{d(F(U_1), F(U_2))}{d(U_1, U_2)}. 
\end{equation}
  By   Condition~(\ref{pfant}) 
   there exist a sequence of fibers  $U_{\lambda_i}\not=U_1$  parallel to $U_1$ such that 
$U_{\lambda_i}$  converges to $U_1$.    For   
sufficiently large $i$,
    we have   
$$\min\{d(p,q), d(p_0, q_0)\}>L\cdot d(U_1, U_{\lambda_i})$$ and 
$$\min\{d(F(p), F(q)), d(F(p_0), F(q_0))\}>L\cdot  d(F(U_1), F(U_{\lambda_i})).$$
 Now Lemma~\ref{l1} applied to $U_1, U_{\lambda_i}$ and $p, q$ 
   yields   
$$\frac{1}{K_1}  \cdot \frac{d(F(U_1), F(U_{\lambda_i}))}{d(U_1, U_{\lambda_i})}\le  \frac{d(F(p), F(q))}{d(p, q)}
\le  K_1\cdot   \frac{d(F(U_1), F(U_{\lambda_i}))}{d(U_1, U_{\lambda_i}))}. $$
 Similarly,  we have 
$$\frac{1}{K_1}  \cdot \frac{d(F(U_1), F(U_{\lambda_i}))}{d(U_1, U_{\lambda_i})}\le  \frac{d(F(p_0), F(q_0))}{d(p_0, q_0)}
\le  K_1\cdot   \frac{d(F(U_1), F(U_{\lambda_i}))}{d(U_1, U_{\lambda_i}))}. $$
 It follows that 
\b{equation}\label{e6}
\frac{1}{K_1^2}\cdot    \frac{d(F(p_0), F(q_0))}{d(p_0, q_0)} \le \frac{d(F(p), F(q))}{d(p, q)}
\le  K_1^2 \cdot  \frac{d(F(p_0), F(q_0))}{d(p_0, q_0)}.
\end{equation}
  Now (\ref{e4})  follows from (\ref{e5})  and (\ref{e6}).

\end{proof}

Lemma~\ref{l2} says that for any fiber $U$  of $X$, the restriction 
 $F|_{U}$ is a $(K_1^3, C)$-quasi-similarity
  for some constant $C>0$ that  may depend on the fiber $U$. 
  Next lemma states that  this constant $C$ can be chosen to be independent of the fiber.

\b{Le}\label{l3}
 There is some constant $C>0$ such that 
$F|_{U}$ is a $(K_2, C)$-quasi-similarity  for any fiber $U$ of $X$, where $K_2$ 
 is a  constant depending only on  $\eta$, $\alpha$,   $L$  and the quasimetric constants of $X$, $Y$. 
\end{Le}

\b{proof}
Fix  a fiber $U_0$ (in $X$)  and let $U$ be an arbitrary fiber in $X$. 
 By Lemma~\ref{l2}, there are constants $C_0, C>0$ such that 
$F|_{U_0}$ is a $(K_1^3, C_0)$-quasi-similarity  and 
$F|_{U}$ is a $(K_1^3, C)$-quasi-similarity. It suffices to show that there is a constant $D$ depending only on  $\eta$, $\alpha$,
     $L$  and the quasimetric constants of $X$, $Y$
 such that $C_0/D\le C\le D C_0$.

  Let $M\ge 1$ be a quasimetric constant for both $X$ and $Y$. 
 Fix $x\in U$ and $x_0\in U_0$,    and pick  $y\in U$ and
 $y_0\in U_0$ such that 
$$d(x_0, y_0)=d(x,y)=10M\cdot d(x, x_0).$$
  The generalized triangle inequality in $X$ applied to 
 $x, x_0, y_0$ implies
$$\frac{1}{2M}\cdot d(x, y)=\frac{1}{2M}\cdot d(x_0, y_0)\le d(x, y_0)\le 2M\cdot d(x_0, y_0)=2M\cdot d(x, y).$$
 The quasisymmetry condition now implies  
$$d(F(x_0), F(y_0))/{\eta(2M)}\le d(F(x), F(y_0))\le \eta(2M)\cdot d(F(x_0),  F(y_0))$$ and 
$$d(F(x),  F(y))/{\eta(2M)}\le d(F(x), F(y_0))\le \eta(2M) \cdot  d(F(x), F(y)).$$
  It follows that 
$$\frac{1}{(\eta(2M))^2}\cdot d(F(x_0), F(y_0))\le d(F(x), F(y))\le (\eta(2M))^2 \cdot d(F(x_0),  F(y_0)).$$
This together with the quasi-similarity property of 
$F|_{U_0}$ and $F|_{U}$  implies 
 $C_0/D\le C\le D C_0$, where $D=(\eta(2M))^2K_1^6$.  
\end{proof}

\b{Le}\label{l4}
$F$ is a $(K, C)$-quasi-similarity, where $K$ depends only on  $\eta$,  $\alpha$,
     $L$  and the quasimetric constants of $X$, $Y$. .
\end{Le}

\b{proof}
 Let $p, q\in X$ be arbitrary. 
  Let $U, V$ be fibers in $X$ such that $p\in U$, $q\in 
 V$ (we may have $U=V$).  
    Pick $x\in U$ such that $d(p, x)=d(p, q)$.  
   Then 
  $$d(F(p), F(q))\le \eta(1)\cdot  d(F(p), F(x))\le \eta(1) K_2 C\cdot  d(p, x)=\eta(1)K_2  C\cdot d(p,q).$$
  So we have the upper bound for $d(F(p), F(q))$.    The same  argument applied to $F^{-1}$   yields 
    the lower
  bound for $d(F(p), F(q))$.   
\end{proof}

The proof of Theorem~\ref{main1} is now complete. 

\vspace{5mm}

\noindent
{\bf{An  Example}}

Recall that  two closed subsets $U$ and $V$ of a quasimetric space  are defined to be parallel if
 $d(u, V)=d(v, U)$ for any $u\in U, v\in V$. 
The following example shows that the conclusion of 
Theorem~\ref{main1}   fails if we replace  \lq\lq parallel" with the weaker condition
 $HD(U, V)<\infty$.

Let $X=Y=\C$ be the  complex  plane with the usual metric. 
 It is well known that for any $\alpha>-1$, the map
 $f_\alpha: \C\ra \C$ given by   $f_\alpha(z)=|z|^{\alpha}z$
 is a quasisymmetric   map.    Let $\alpha>-1$,  $\alpha\not=0$, define a map 
 $F_\alpha: \C\ra \C$ as follows:
 \begin{equation*}
F_\alpha(z)= \left\{
\begin{array}{rl}
f_\alpha(z) &   \text{if} \;  |z|\le 1\\
z  &  \text{if} \; |z|\ge 1.
\end{array}  \right.  
\end{equation*}
 Then it is clear that $F_\alpha$ is a quasiconformal map, and hence is quasisymmetric.
   The fibers in $X$ are horizontal lines, and the fibers in $Y$ are images of horizontal lines  under $F_\alpha$.   
  A  direct   calculation  shows   that all fibers of $Y$ are $10$-biLipschitz  to the real line.  So  Condition (\ref{fasfetugs})    is satisfied by both $X$ and $Y$.  
 All  the conditions in 
 Definition \ref{fibered} 
    are satisfied, provided we replace
  \lq\lq parallel"   with the weaker condition  $HD(U, V)<\infty$. 
  However,  $F_\alpha$ is not biLipschitz due to the distorsion around the origin.

\section{Applications}\label{app}

\subsection{Quasisymmetric maps on Carnot groups with reducible first stratum}\label{liouville}

In this section we  use Theorem~\ref{main1}  to prove Theorem~\ref{main2}.    Specifically,
    we 
show that every quasisymmetric map on Carnot groups with reducible first stratum  
is   biLipschitz. 

Let   $G$ be a Carnot group with Lie algebra  
$\mathfrak g=V_1\oplus \cdots \oplus V_s$.
Assume there is 
    a non-trivial proper linear subspace   
$W_1\subset  V_1$ that is invariant under the action of the
group of strata-preserving automorphisms of $\mathfrak g $.
   Let  $\mathfrak h$ be the subalgebra of $\mathfrak{g}$ generated by $W_1$.  
We denote the connected   Lie  subgroup of $G$   with Lie algebra  $\mathfrak h$ 
by $H$ and refer to it as   the {\em subgroup generated} by $W_1$.
 Notice that $H$ is also a Carnot group and its lie algebra can be written as 
   $\mathfrak h=W_1\oplus \cdots \oplus  W_s$.    In general, there is some 
  integer  
 $1\le \bar s\le  s$ such that  $W_{\bar s}\not=0$ and $W_j=0$ for $j> \bar s$.  
  Let  $\langle\cdot,\cdot\rangle$    be   an inner product   
  on $V_1$  and $d$    the   left-invariant   sub-Riemannian   Carnot  metric   on $G$     determined by   $\langle\cdot,\cdot\rangle$.


We recall the following: 

\b{Prop}\label{p3}   \e{(Proposition 3.4, \cite{X2})}
 Let  $G$  and $G'$ be two Carnot groups.
%
%
  Let $W_1 $  and $ W'_1$ be two subspaces of the first strata of the stratifications of  ${\rm Lie}(G)$ and
    ${\rm Lie}(G')$,
respectively.
Let $H$ and $H'$ be the groups generated by $W_1$ and $ W'_1$, respectively.
   Let  
  $F:   G \ra G'$ be a quasisymmetric     map.  
  If     $dF(x)(W_1)\subset W'_1$ for a.e. $x\in  G$,   
     then   $F$ sends  each left coset $U$ of $H$ 
  into a  left  coset    of $H'$.  

  \end{Prop}

Now let $F: G \ra G$ be an $\eta$-quasisymmetric map. 
 By Pansu's differentiability theorem,  at  almost everywhere  $x\in G$ the map  $F$ is Pansu differentiable   and   
 the Pansu differential 
 $dF(x): \mathfrak g\ra \mathfrak g$ is a strata-preserving automorphism.  
 The assumption on $W_1$ 
 in Theorem~\ref{main2}  
    implies that 
 $dF(x)(W_1)=W_1$.  Now  Proposition~\ref{p3} implies that $F$ sends left cosets of $H$ to left cosets of $H$, where $H$ is the subgroup generated by $W_1$.


Theorem \ref{main2} shall follow from Theorem \ref{main1} once we verify the conditions in Definition~\ref{fibered}.  Here the fibers are left cosets of $H$.
Since $H$ is also a Carnot group,  $(H, d)$ is 
  a   geodesic space.
  So Condition \eqref{fasfetugs}   is satisfied for $\alpha=L=1$.    Conditions    
\eqref{pfant}, \eqref{fhpd}  and   \eqref{npfd}   will be verified in Lemmas 
\ref{nonisoparallel2},      \ref{nonzerod}   and \ref{nonparalleldiverge}, respectively.





 Given a Lie algebra 
$\mathfrak n$     and a subalgebra 
$\mathfrak{s}\subset \mathfrak n$,  
  the normalizer of $\mathfrak{s}$ in $\mathfrak n$ is  defined by:
 $$ \mathcal N_{\mathfrak n}(\mathfrak{s})=\{X\in \mathfrak n: [X, \mathfrak{s}]\subset \mathfrak{s}\}.$$
 It is easy to see that $ \mathcal N_{\mathfrak n}(\mathfrak{s})$ is a Lie subalgebra of  $\mathfrak n$ and 
$\mathfrak{s}$ is an ideal in $ \mathcal N_{\mathfrak n}(\mathfrak{s})$.  
 Let $N$ be a connected Lie group with Lie algebra $\mathfrak n$.
 Let $S$ be the Lie subgroup of $N$ with Lie algebra $\mathfrak s$. 
 Let $d$ be a left invariant distance on $N$. 
  Lemmas \ref{nonisoparallel2}  and \ref{nonzerod}  are valid  for  all connected and simply connected nilpotent Lie groups  $N$  with a left invariant distance. 
 In particular, they hold true for the Carnot group $G$ and left cosets of $H$.  

\b{Le}\label{nonisoparallel2}
For any  proper   Lie  subalgebra $\mathfrak{s}$ of a  nilpotent Lie algebra
$\mathfrak n$, we have
\begin{equation}\label{nonisoparallel}
 \mathcal N_{\mathfrak n}(\mathfrak{s})\not=\mathfrak{s}.
 \end{equation}
Consequently, 
for any left coset $U$ of $S$, there exist a sequence of left cosets $U_i$ of $S$
that are  parallel to $U$ and converge to $U$. 
\end{Le}

\b{proof} For the first claim 
   we   induct on the degree of nilpotency. The    claim holds  when $\mathfrak  n $ is abelian 
 since in this case  $ \mathcal N_{\mathfrak n}(\mathfrak{s})=\mathfrak n$.
Suppose the   claim  holds for all  $k$-step nilpotent Lie algebras.    Let 
$\mathfrak n$ be  $(k+1)$-step.   We assume  $\mathcal N_{\mathfrak n}(\mathfrak{s})=\mathfrak{s}$ and will derive a contradiction from this.
    Notice the center $C(\mathfrak n)$ of $\mathfrak n$ lies in
      $\mathcal N_{\mathfrak n}(\mathfrak{s})=\mathfrak{s}$. 
 So we have a proper   Lie  subalgebra $\mathfrak{s}/C(\mathfrak n)$ of
 the $k$-step nilpotent Lie algebra $\mathfrak{n}/C(\mathfrak n)$.  By the induction hypothesis,
 we have  $  \mathcal N_{\mathfrak n/C(\mathfrak n)}(\mathfrak{s}/C(\mathfrak n))\not=\mathfrak{s}/C(\mathfrak n)$.
Let $\pi:  \mathfrak n\ra {\mathfrak n}/C(\mathfrak n)$  be the natural projection. 
 It is easy to check that 
$\pi^{-1}(  \mathcal N_{\mathfrak n/C(\mathfrak n)}(\mathfrak{s}/C(\mathfrak n)))= \mathcal N_{\mathfrak n}(\mathfrak{s})$. 
 Since $\pi^{-1}(\mathfrak{s}/C(\mathfrak n))=\mathfrak{s}$,  
  we obtain  $\mathcal N_{\mathfrak n}(\mathfrak{s})\not=\mathfrak{s}$,
  contradicting the assumption  $\mathcal N_{\mathfrak n}(\mathfrak{s})=\mathfrak{s}$.

  Regarding the second part of the lemma, let  $K$   be   the   connected subgroup of $N$ with Lie algebra    
$ \mathcal N_{\mathfrak n}(\mathfrak s)$.   Since $\mathfrak s$ is an ideal in $ \mathcal N_{\mathfrak n}(\mathfrak s)$,
   we see that $S$ is a  normal  subgroup of   $K$.  
From  \eqref{nonisoparallel} we deduce that $K$ contains $S$ properly. 
There exist $k_i\in K\backslash S$ such that  $k_i\ra e$. 
 So for any left coset  $U=gS$, the  left cosets $U_i=gk_iS$  
  converge to $U$.  Finally we claim that  $U_i$ is parallel to $U$. Indeed,   
  for any $p=gs_0\in U$,  
$$d(p, U_i)=d(g s_0,  gk_i S)=d(e, s_0^{-1}k_iS)
=d(e, k_i\cdot k_i^{-1}s_0^{-1}k_iS)=d(e, k_i S);$$
  and for any $q=gk_i s_0\in U_i$,
  \begin{align*} d(q, U) & =d(gk_i s_0,  gS)\\
& =d(Sk_is_0, e)=d(Sk_is_0k_i^{-1}\cdot k_i, e)=d(Sk_i,  e)=d(e,  Sk_i)=d(e, k_iS),
\end{align*}
  where in the last equality we used the fact that  $Sk_i=k_iS$.
\end{proof}

Notice that the last part of the preceding proof shows that for any $g\in N$ and any $k\in K$, the two left cosets $gS$ and $gkS$ are parallel.

Next lemma says that two different left costs of $S$ can never  get arbitrarily close. 

\b{Le}\label{nonzerod}
If $U$ and $V$ are two distinct left cosets of $S$,   then $d(U,V )>0$.

\end{Le}

\b{proof}  Set  $\mathfrak{s}_0=\mathfrak s$   
and define inductively
 $\mathfrak{s}_j= \mathcal N_{\mathfrak n}(\mathfrak{s}_{j-1})$. 
 From    \eqref{nonisoparallel} we have that 
   $\mathfrak{s}_j$ properly contains  $\mathfrak{s}_{j-1}$  unless $\mathfrak{s}_{j-1}=\mathfrak n$. 
 Since $\mathfrak n$ is finite dimensional, there is some $k$ such that $\mathfrak{s}_{k-1}\not=\mathfrak n$
  and  
 $\mathfrak{s}_k=\mathfrak n$. Let $S_j$ be the   connected Lie  subgroup
 of $N$ with Lie algebra 
$\mathfrak{s}_j$. Then $S_{j-1}$  is a proper normal sugroup of  $S_j$  for $j\le k$.
  It follows that for any $g\in S_j\backslash\ S_{j-1}$, the two left cosets 
 $S_{j-1}$  and $gS_{j-1}$  are parallel. Hence
 $d(S_{j-1}, gS_{j-1})>0$.

Now let $U$ and $V$ be   two distinct left cosets of $S$. After applying a left translation we may assume $U=S$ and $V=gS$ for some $g\notin S$.
There is some $j\le k$ such that $g\in S_j\backslash  S_{j-1}$.
It follows that 
$$d(U,V)=d(S, g S)\ge d(S_{j-1}, gS_{j-1})>0,$$
since
$S\subset S_{j-1}$.
\end{proof}

  For the proof of the next lemma we will work under the assumption of Theorem \ref{main2}.

Let $\langle\cdot,\cdot\rangle$ be an inner product on $\mathfrak g=V_1\oplus \cdots \oplus V_s$ such that $V_i$ and $V_j$ are perpendicular to each other for $i\not=j$.   
 Define a \lq\lq quasi norm" on $\mathfrak g$ as follows.
 $$||\sum_i x_i||=\sum_i \langle x_i, x_i\rangle ^{\frac{1}{2i}},  \; \text{where}\; x_i\in V_i.$$
   Next define  $\bar d(x,y)=||(-x)*y||$ for $x, y\in \mathfrak g$. 
 It is well known that every Carnot metric $d$  on $\mathfrak g$  is  
  biLipschitz equivalent  to $\bar d$:        there exists some constant $C_0\ge 1$ such that
 $$\frac{1}{C_0} \cdot \bar d(x,y)\le d(x,y)\le C_0 \cdot \bar d(x,y)$$
  for   any  $x, y\in \mathfrak g$. 
  Here we identified $G$ and $\mathfrak g$ via the exponential map.

\b{Le}\label{nonparalleldiverge}
If $U$ and $V$ are two  left cosets of $H$ that are not parallel, then $HD(U, V)=\infty$.
\end{Le}

\b{proof}  Suppose $U$ and $V$ are not parallel and $HD(U, V)<\infty$.
 We shall obtain a contradiction. 
After applying a left translation we may assume $U=H$ and $V=g H$ for some 
 $g\in G\backslash H$.     The last part of the proof of Lemma \ref{nonisoparallel2} shows that 
 for any $k$ in the normalizer $K$  of $H$  the two fibers
 $H$ and $kH=Hk$ are parallel. So    $g\notin K$.
Below we will  identify $\mathfrak g$ with $G$  via the exponential map
   and do calculations  in the Lie algebra using BCH formula.  
 So we let  $Y\notin \mathcal{N}_\mathfrak{g}(\mathfrak h)$ be such that 
 $HD(\mathfrak h, Y*\mathfrak h)<\infty$.  Since 
  $$HD(Y*\mathfrak h,   \;  Y*\mathfrak h *(-Y))\le d(e, -Y),$$
 we see that $HD(\mathfrak h, \; Y*\mathfrak h *(-Y))\le C$  for some constant $C>0$.

We write $Y=\sum_{i=1}^s Y_i$  with $Y_i\in V_i$  and let $j\ge 1$ be the index such that 
  $Y_j\notin \mathcal{N}_\mathfrak{g}(\mathfrak h)$  and 
$Y_i\in \mathcal{N}_\mathfrak{g}(\mathfrak h)$ for all  $i<j$. 
 By replacing  $Y$ with $Y*(-Y_1-\ldots-Y_{j-1})$, we may assume $Y_i=0$ for $ i<  j$.     Since  $[Y_j, \mathfrak h]\not\subset \mathfrak h$  and 
$\mathfrak h$ is generated by the first  stratum  $W_1$,
  there is some $X\in W_1$ such that $[Y_j, X]\notin \mathfrak h$. 
  Notice that $[Y_j, X]\in V_{j+1}$.  
Since 
$HD(\mathfrak h,  Y*\mathfrak h *(-Y))\le C$,
  for any $t\in \R$, there is some $X'\in \mathfrak h$ ($X'$ may depend on $t$) such that
 \begin{equation}\label{e40}
d((-X')*Y*(tX)*(-Y), 0) =    d(Y*(tX)*(-Y), X')\le C. 
\end{equation}

We next calculate $Z:=Y*(tX)*(-Y)$ and $A:=(-X')*Z$. 
We have $$Z=e^{{\rm ad}\, Y} (tX)=t\{X+[Y, X]+\cdots +\frac{1}{k!}({\rm ad} \,Y)^k (X)+\cdots+\frac{1}{(s-1)!}({\rm ad}\, Y)^{s-1} (X)\},$$
 where ${\rm ad}\, Y: \mathfrak g\ra \mathfrak g$ is  the linear map given by ${\rm ad}\, Y (B)=[Y, B]$ for $B\in \mathfrak g$.
  Write $Z=Z_1+\cdots + Z_s$ with $Z_i\in V_i$.  
Since  $Y_i=0$ for $ i<j$,  we have 
$$\frac{1}{k!}({\rm ad} \,Y)^k (X)\in V_{j+2}\oplus \cdots \oplus V_s\; \; \text{for  all}\;\; k\ge 2.$$
 So the terms $Z_i$ with $i\le j+1$ is determined by $t(X+[Y, X])$. 
  Since 
$X\in V_1$  
we have
 $Z_1=tX$, $Z_i=0$ for $2\le i\le j$ and $Z_{j+1}=t[Y_j, X]$.

 Write $A=A_1+\cdots + A_s$ with $A_i\in V_i$.  
  Notice that  $[-X', Z]$ and all the iterated brackets in the BCH formula for 
 $(-X')*Z$ are the sum of an element of $\mathfrak h$ and an element in 
 $V_{j+2}\oplus \cdots \oplus V_s$.   
  So $A_{j+1}$ is completely determined by $-X'+Z$. 
It follows that 
 $A_{j+1}$ is the sum of an   element of 
 $W_{j+1}$ and  $t[Y_j, X]$ (recall $\mathfrak h=W_1\oplus \cdots \oplus W_s$).  
Write $[Y_j, X]=B+B^\perp$, where $B\in W_{j+1}$ and $B^\perp\in V_{j+1}$ is perpendicular to $W_{j+1}$ with respect to the inner product $<, >$ on $\mathfrak g$.
   Since $[Y_j, X]\notin \mathfrak h$, we have 
$B^\perp\not=0$.  
Hence 
$A_{j+1}$ equals the sum of $tB^\perp$ and an element of $W_{j+1}$.
   It follows that 
 $$d((-X')*Y*(tX)*(-Y), 0) \ge \frac{1}{C_0}\cdot  ||(-X')*Y*(tX)*(-Y)||\ge \frac{1}{C_0}\cdot |tB^\perp|^{\frac{1}{j+1}}\ra \infty,$$
 as $t\ra \infty$ since $B^\perp\not=0$.  
  This contradicts  (\ref{e40}).
\end{proof}

We have verified all  the conditions in Definition \ref{fibered}  and 
the proof of Theorem \ref{main2} is now complete.

\begin{remark}\label{remark4}
Lemma \ref{nonparalleldiverge}   is equivalent  to the following  statement:   if the Hausdorff  distance between $H$ and $gHg^{-1}$ is finite, then $g$ lies in the normalizer of $H$.
  We believe  this  is true for any  left invariant distance on any connected, simply connected nilpotent Lie group $N$ and any proper Lie subgroup $H$.  
\end{remark}

\subsection{Other applications}\label{other}

In this subsection we give a couple of other applications of Theorem \ref{main1}.
 The first is to quasiconformal maps of Heisenberg groups that send 
    vertical lines to vertical lines. The second is to quasisymmetric maps of 
 ideal boundary of certain  amenable hyperbolic locally compact groups. 

 We first show that quasiconformal maps of the Heisenebrg groups 
 $H^n=\R^{2n}\times \R$   that  permute vertical lines are biLipschitz. 
 These maps are lifts of those biLipschitz maps of $\R^{2n}$ that preserve the standard sympletic form on $\R^{2n}$.  
 The reader is referrd to \cite{T} and \cite{BHT}   for   more details.

\b{Prop}\label{hei}
Let $F: H^n\ra H^n$ be a quasiconformal map of the Heisenberg group. If 
 $F$ maps vertical lines to vertical lines, then $F$ is biLipschitz.

\end{Prop}

\b{proof}  Recall that by  \cite[Theorem 4.7]{HK}, $F$ is quasisymmetric.
The fibers are the vertical lines.   There is some $L\ge 1$ such that the vertical lines 
 are  $L$-biLipschitz to 
 $(\R, |\cdot|^{\frac{1}{2}})$, where $|\cdot|$ is the usual metric on 
 $\R$.  So Condition  (\ref{fasfetugs})    in Definition~\ref{fibered}  is satisfied for $\alpha=1/2$ and $L$. It is easy to see
  that   the other three conditions 
   in   Definition~\ref{fibered}   are also satisfied.     By Theorem~\ref{main1}, 
  $F$ is biLipschitz. 
\end{proof}

We next give a simple proof of  Dymarz's theorem.

\begin{theorem}[{\cite[Theorem 1]{D}}]
 Let $N$ be a Carnot group equipped with a left invariant Carnot metric $d$
  and $\Q_m$  ($m \ge 2$) the $m$-adics  with standard metric $d_m$. Then 
 every quasisymmetric map $N\times \Q_m \ra N\times \Q_m$ is biLipschitz.
  Here the metric on $N\times \Q_m$ is   given by 
 $d'((x_1, y_1), (x_2, y_2))=\max \{d(x_1, x_2), d_m(y_1, y_2)\}$.

\end{theorem}

Recall  that the standard metric $d_m$ on $\Q_m$ is given by:
  $$d_m(\sum a_i m^i, \sum b_i m^i)= m^{-(k+1)},  $$
where $k$ is the smallest    index for which $a_i \neq b_i.$

\begin{proof}
  The fibers of $N\times \Q_m$ are the subsets  $N\times \{y\}$ ($y\in \Q_m$).
Notice that $(\Q_m, d_m)$ is perfect and totally disconnected.  So the fibers of 
 $N\times \Q_m$ are
 exactly the connected components of  $N\times \Q_m$.
  Hence every quasisymmetric map $F$ of 
$N\times \Q_m$    permutes the  fibers.  
 Each fiber is  isometric to a Carnot group and so is geodesic.  Hence  Condition 
(\ref{fasfetugs})   is satisfied for $\alpha=L=1$.
 The perfectness of $\Q_m$ implies that  parallel  fibers are not isolated (Condition (\ref{pfant})).  
   Distance between distinct fibers is clearly  positive  (Condition (\ref{fhpd})).
All fibers are parallel, so Condition  (\ref{npfd})  is vacuous. 
By Theorem \ref{main1} every  quasisymmetric map 
 of $N\times \Q_m$  is biLipschitz. 
\end{proof}

 \addcontentsline{toc}{subsection}{References}

\noindent Addresses:

\noindent  Enrico Le Donne:  Department of Mathematics and Statistics,  
P.O. Box 35 (MaD), 
FI-40014 University of Jyvaskyla,
Finland.  \hskip .4cm E-mail:   ledonne@msri.org

\vspace{5mm}

\noindent Xiangdong Xie: Dept. of Mathematics  and   Statistics,   Bowling Green  State  University, 
  Bowling Green,  OH,   U.S.A.\hskip .4cm E-mail:   xiex@bgsu.edu

\end{document}